\documentclass{amsart}

\newtheorem*{theorem*}{Theorem}
\newtheorem{theorem}{Theorem}
\newtheorem{thmx}{Theorem}

\newtheorem{lem}[thmx]{Lemma}

\newtheorem{remark}{Remark}

\begin{document}
\title[$PGT$ on $PSL(2,\mathbb{Z})$]{$PGT$ on $PSL(2,\mathbb{Z})$. A short
proof}
\author{Muharem Avdispahi\'{c}}
\address{University of Sarajevo, Department of Mathematics, Zmaja od Bosne
33-35, 71000 Sarajevo, Bosnia and Herzegovina}
\email{mavdispa@pmf.unsa.ba}

\maketitle

\begin{abstract}
Taking the Iwaniec explicit formula as a starting point, we give a
short proof of a more precise $\frac{2}{3}$ bound for the exponent in
the error term of the Gallagher-type prime geodesic theorem for the
modular surface.
\end{abstract}

\section{Introduction}

Let $\Gamma = PSL\left( 2,%
\mathbb{Z}
\right) $ be the modular group and $Z_{\Gamma }$ its Selberg function
defined by
\begin{equation*}
Z_{\Gamma }(s)=\underset{\{P_{0}\}}{\prod }\underset{k=0}{\overset{\infty }{%
\prod }}(1-N(P_{0})^{-s-k})\text{, }\text{Re}(s)>1\text{,}
\end{equation*}
and meromorphicaly continued to the whole complex plane. The product is over
hyperbolic conjugacy classes in $\Gamma$. The primitive conjugacy classes $%
P_0$ correspond to primitive closed geodesics on the modular surface $\Gamma
\setminus \mathcal{H}$, where $\mathcal{H}$ is the upper half-plane equipped
with the hyperbolic metric. The length of the primitive closed geodesic
joining two fixed points, necessarily the same for all representatives of a
class $P_{0}$, equals $\log (N(P_0))$. We are interested in distribution of
these geodesics, i.e., in the number $\pi _{\Gamma}(x)$ of classes $P_{0}$
such that $N(P_{0})\leq x$, for $x>0$.

It is believed that the error term in the prime geodesic theorem
\begin{equation*}
\pi_{\Gamma }\left( x\right) \approx \int_{0}^{x}\frac{dt}{\log t} \ \
\left( x\rightarrow \infty \right)
\end{equation*}
 is $O(x^{\frac{1}{2}+\varepsilon})$ since $Z_{\Gamma }$ satisfies the Riemann hypothesis.

Namely, there exists a broad analogy between the prime geodesic theorem and the prime number theorem based on the role played by the distribution of
zeros of the Selberg zeta function resp. Riemann zeta function in two respective contexts.

However, due to the fact that $Z_{\Gamma }$ is a meromorphic function of
order $2$, as opposed to the Riemann zeta which is of order $1$, the best
presently available estimate for the exponent in the error term is $\frac{25}{36}+\varepsilon$, obtained by
Soundararajan and Young \cite{SY:2013}. The generalized Lindel\"{o}f
hypothesis for Dirichlet $L$-functions would imply $\frac{2}{3}+\varepsilon $
(see \cite{I1:1984}, \cite{SY:2013}).

Inspired by Gallagher's approach \cite{G1:1980} in the Riemann zeta setting,
we give a short proof of the following theorem.

\begin{theorem}
\label{t1} Let $\Gamma = PSL(2,\mathbb{Z})$ be the modular group and $%
\varepsilon >0$ arbitrarily small. There exists a set $A $ of finite
logarithmic measure such that
\begin{equation*}
\pi _{\Gamma }\left( x\right) = \int_{0}^{x}\frac{dt}{\log t}+O\left( x^{%
\frac{2}{3}}\left( \log \log x\right) ^{\frac{1}{3}+ \varepsilon }\right) \ \
\left( x\rightarrow \infty,\ x\notin A\right) \text{.}
\end{equation*}
\end{theorem}

\section{Preliminaries. The Iwaniec explicit formula and Gallagher's lemma.}

We shall make use of two lemmas that played important role in \cite{I2:1984}
and \cite{G1:1980}. The first one is the Iwaniec explicit formula \cite{I2:1984} with an
error term for the Chebyshev function
\begin{equation*}
\psi _{\Gamma }\left(
x\right) =\underset{N\left( P\right) \leq x}{\sum }\log N\left(
P_{0}\right) =\underset{N\left( P_{0}\right) ^{k}\leq x}{\sum }\log N\left(
P_{0}\right)\text{.}
\end{equation*}

\begin{lem}
\label{lemma1} For $1\leq T\leq \frac{x^{\frac{1}{2}}}{\left(\log x\right)^2}
$, one has
\begin{equation*}
\psi_{\Gamma} \left( x\right) =x+\underset{\left\vert \gamma \right\vert
\leq T}{\sum }\frac{x^{\rho }}{\rho } +O\left(\frac{x}{T} \left(\log
x\right)^2\right)\text{,}
\end{equation*}
where $\rho = \frac{1}{2}+i\gamma$ denote zeros of $Z_{\Gamma}$.
\end{lem}

The second one is Gallagher's lemma \cite{G:1970} that enabled him to reduce
the error term in the prime number theorem under the Riemann hypothesis.

\begin{lem}
\label{lemma2} Let $A$ be a discrete subset of $\mathbb{R}$ and $\theta \in
(0,1)$. For any sequence $c(\nu )\in \mathbb{C}$, $\nu \in A$, let the
series
\begin{equation*}
S\left( u\right) =\underset{\nu \in A}{\sum }c\left( \nu \right) e^{2\pi
i\nu u}
\end{equation*}%
be absolutely convergent. Then%
\begin{equation*}
\int_{-U}^{U}\left\vert S\left( u\right) \right\vert ^{2}du\leq \left( \frac{%
\pi \theta }{\sin \pi \theta }\right) ^{2}\int_{-\infty }^{+\infty
}\left\vert \frac{U}{\theta }\underset{t\leq \nu \leq t+\frac{\theta }{U}}{%
\sum }c\left( \nu \right) \right\vert ^{2}dt\text{.}
\end{equation*}
\end{lem}

\section{ Proof of Theorem}

The goal is to find a proper bound for $\underset{\left\vert \gamma
\right\vert \leq T}{\sum }\frac{x^{\rho }}{\rho }$ in Lemma \ref{lemma1}.
For $x\in \left[ e^{n},e^{n+1}\right) $, we have
\begin{equation}
\overset{e^{n+1}}{\underset{e^{n}}{\int }}\left\vert \underset{\left\vert
\gamma \right\vert \leq T}{{\sum }}\frac{x^{\rho }}{\rho }\right\vert ^{2}dx=%
\overset{e^{n+1}}{\underset{e^{n}}{\int }}x^{2}\left\vert \underset{%
\left\vert \gamma \right\vert \leq T}{{\sum }}\frac{x^{i\gamma }}{\rho }%
\right\vert ^{2}\frac{dx}{x}\ll e^{2n}\overset{e^{n+1}}{\underset{e^{n}}{%
\int }}\left\vert \underset{\left\vert \gamma \right\vert \leq T}{{\sum }}%
\frac{x^{i\gamma }}{\rho }\right\vert ^{2}\frac{dx}{x}\text{.}  \label{eq1}
\end{equation}

Through the substitution $x=e^{n}\cdot e^{2\pi \left( u+\frac{1}{4\pi }\right)
}$, the last integral is transformed into $2\pi \overset{\frac{1}{4\pi }}{\underset{-%
\frac{1}{4\pi }}{\int }}\left\vert \underset{\left\vert \gamma \right\vert
\leq T}{{\sum }}\frac{e^{\left( n+\frac{1}{2}\right) i\gamma }}{\rho }%
e^{2\pi i\gamma u}\right\vert ^{2}du$.

Lemma \ref{lemma2}., with $\theta =U=\frac{1}{4\pi }$ and $c_{\gamma }=\frac{%
e^{\left( n+\frac{1}{2}\right) i\gamma }}{\rho }$ for $\left\vert \gamma
\right\vert \leq T$, $c_{\gamma }=0$ otherwise, implies%
\begin{equation}
\overset{\frac{1}{4\pi }}{\underset{-\frac{1}{4\pi }}{\int }}\left\vert
\underset{\left\vert \gamma \right\vert \leq T}{{\sum }}\frac{e^{\left( n+%
\frac{1}{2}\right) i\gamma }}{\rho }e^{2\pi i\gamma u}\right\vert ^{2}du\leq
\left( \frac{\frac{1}{4}}{\sin \frac{\pi }{4}}\right) ^{2}\overset{+\infty }{%
\underset{-\infty }{\int }}\left( \underset{{\scriptsize
\begin{array}{c}
t<\left\vert \rho \right\vert \leq t+1 \\
\left\vert \gamma \right\vert \leq T%
\end{array}%
}}{{\sum }}\frac{1}{\left\vert \rho \right\vert }\right) ^{2}dt\text{.}
\label{eq2}
\end{equation}

According to the Weyl law, $\underset{{\scriptsize t<\left\vert \rho
\right\vert \leq t+1}}{{\sum }}\frac{1}{\left\vert \rho \right\vert }%
=O\left( 1\right) $. Thus,%
\begin{equation}
\overset{+\infty }{\underset{-\infty }{\int }}\left( \underset{{\scriptsize
\begin{array}{c}
t<\left\vert \rho \right\vert \leq t+1 \\
\left\vert \gamma \right\vert \leq T%
\end{array}%
}}{{\sum }}\frac{1}{\left\vert \rho \right\vert }\right) ^{2}dt=O\left(
\overset{T+1}{\underset{\frac{1}{2}}{\int }}dt\right) =O\left( T\right)
\text{.}  \label{eq3}
\end{equation}

Combining \eqref{eq1}, \eqref{eq2} and \eqref{eq3}, we get%
\begin{equation}
\overset{e^{n+1}}{\underset{e^{n}}{\int }}\left\vert \underset{\left\vert
\gamma \right\vert \leq T}{{\sum }}\frac{x^{\rho }}{\rho }\right\vert
^{2}dx=O\left( e^{2n}T\right) \text{.}  \label{eq4}
\end{equation}

Let $A_{n}=\left\{ x\in \left[ e^{n},e^{n+1}\right) :\left\vert \underset{%
\left\vert \gamma \right\vert \leq T}{{\sum }}\frac{x^{\rho }}{\rho }%
\right\vert >x^{\frac{1}{2}}T^{\frac{1}{2}}\left( \log x\right) ^{\frac{1}{2}%
}\left( \log \log x\right) ^{\frac{1}{2}+\frac{3}{2}\varepsilon }\right\} $. Its \\
\noindent logarithmic measure is controlled by

\begin{eqnarray*}
\mu ^{\ast }A_{n} &=&\underset{A_{n}}{\int }\frac{dx}{x}=\underset{A_{n}}{%
\int }xT\left( \log x\right) \left( \log \log x\right) ^{1+3\varepsilon }%
\frac{dx}{x^{2}T\left( \log x\right) \left( \log \log x\right)
^{1+3\varepsilon }} \\
&\leq &\frac{1}{e^{2n}Tn\left( \log n\right) ^{1+3\varepsilon }}\overset{%
e^{n+1}}{\underset{e^{n}}{\int }}\left\vert \underset{\left\vert \gamma
\right\vert \leq T}{{\sum }}\frac{x^{\rho }}{\rho }\right\vert ^{2}dx\ll
\frac{e^{2n}T}{e^{2n}Tn\left( \log n\right) ^{1+3\varepsilon }}=\frac{1}{%
n\left( \log n\right) ^{1+3\varepsilon }}\text{.}
\end{eqnarray*}%
Thus,%
\begin{equation}
\left\vert \underset{\left\vert \gamma \right\vert \leq T}{{\sum }}\frac{%
x^{\rho }}{\rho }\right\vert \leq x^{\frac{1}{2}}T^{\frac{1}{2}}\left( \log
x\right) ^{\frac{1}{2}}\left( \log \log x\right) ^{\frac{1}{2}+\frac{3}{2}%
\varepsilon }  \label{eq5}
\end{equation}%
outside a set $A=\cup A_{n}$ of finite logarithmic measure.

The optimal choice is $T\approx \frac{x^{\frac{1}{3}}\log x}{\left( \log
\log x\right) ^{\frac{1}{3}+\varepsilon }}$. Then, Lemma \ref{lemma1}. and
the relation \eqref{eq5} yield%
\begin{equation}
\psi _{\Gamma }\left( x\right) =x+O\left( x^{\frac{2}{3}}\log x\left( \log
\log x\right) ^{\frac{1}{3}+\varepsilon }\right) \text{ \ \ }\left(
x\rightarrow \infty \text{, }x\notin A\right) \text{.}  \label{eq6}
\end{equation}%
From \eqref{eq6}, we obtain the assertion of Theorem \ref{t1} in a standard way,
making use of the expressions

\begin{equation*}
\pi _{\Gamma }\left( x\right) =\int_{2}^{x}\frac{1}{\log t}d\theta _{\Gamma
}\left( t\right) \text{ and }\psi _{\Gamma }\left( x\right) =\overset{\infty
}{\underset{n=1}{\sum }}\theta _{\Gamma }\left( x^{\frac{1}{n}}\right) \text{%
,}
\end{equation*}%
where $\ \theta _{\Gamma }\left( x\right) =\underset{N\left( P_{0}\right)
\leq x}{\sum }\log N\left( P_{0}\right) $.

\begin{remark}
In the case of cofinite Fuchsian groups  $\Gamma \subset PSL(2,\mathbb{%
\mathbb{R}
})$, the best unconditional estimate of the remainder in the prime geodesic
theorem is still Randol's $O\left( \frac{x^{\frac{3}{4}}}{\log x}\right) $
\cite{R:1977}. (See \cite{A:2017} for a passage from Hejhal's proof \cite[Th. 6.19]{H:1976} of Huber's \cite{H1:1961, H2:1961} bound $O\left( \frac{x^{\frac{3}{%
4}}}{(\log x)^{\frac{1}{2}}}\right) $ to Randol.) Its analogue is also
valid in higher dimensions \cite{AG:2012}. Randol's estimate can be reduced
to $\frac{7}{10}+\varepsilon $ outside a set of finite logarithmic measure
\cite{A1:2017}, what coincides with the Luo-Sarnak \cite{LS:1995} bound that unconditionally
holds for $\Gamma =PSL(2,\mathbb{%
\mathbb{Z}
})$.
\end{remark}

\begin{remark}
Using a more involved method of Soundararajan and Young, we were able to
prove \cite{A2:2017} that under the assumption of the generalized Lindel{\"o}f
hypothesis one has
\begin{equation*}
\pi _{\Gamma }\left( x\right) =\int_{0}^{x}\frac{dt}{\log t}+O\left( x^{\frac{5}{8}%
+\varepsilon }\right)
\end{equation*}%
as $x\rightarrow \infty$ outside a set of finite logarithmic measure.
\end{remark}


\begin{thebibliography}{0}
\bibitem{A:2017} M. Avdispahi\'{c}, On Koyama's refinement
of the prime geodesic theorem, arXiv:1701.01642.

\bibitem{A1:2017} M. Avdispahi\'{c}, Prime geodesic theorem of Gallagher type, arXiv:1701.02115.

\bibitem{A2:2017} M. Avdispahi\'{c}, Prime geodesic theorem for the modular surface, arXiv:1702.01699.

\bibitem{AG:2012} M. Avdispahi\'{c} and D\v{z}. Gu\v{s}i\'{c}, On the error term in the prime geodesic theorem, \textit{Bull. Korean Math. Soc.} \textbf{49} (2012), no. 2, 367--372.

\bibitem{G:1970} P. X. Gallagher, A large sieve density estimate near $\sigma =1$, \textit{Invent. Math.} \textbf{11} (1970), 329--339.

\bibitem{G1:1980} P. X. Gallagher, Some consequences of the Riemann hypothesis, \textit{Acta Arith.} \textbf{37} (1980), 339--343.

\bibitem{H:1976} D. A. Hejhal, \textit{The Selberg trace formula for $%
PSL(2,R)$. Vol I}, Lecture Notes in Mathematics, Vol \textbf{548}. Springer, Berlin 1976.

\bibitem{H1:1961} H. Huber, Zur analytischen Theorie hyperbolischer
Raumformen und Bewegungsgruppen II, \textit{Math. Ann.} \textbf{142} (1961), 385-398.

\bibitem{H2:1961} H. Huber, Nachtrag zu \cite{H1:1961}, \textit{Math. Ann.} \textbf{143}
(1961), 463-464.

\bibitem{I1:1984} H. Iwaniec, Non-holomorphic modular forms and their applications. In \textit{Modular forms} (Durham, 1983), 157--196, Ellis Horwood Ser. Math. Appl.: Statist. Oper. Res., Horwood, Chichester 1984.

\bibitem{I2:1984} H. Iwaniec, Prime geodesic theorem, \textit{J. Reine Angew. Math.} \textbf{349} (1984), 136--159.

\bibitem{K:2016} S. Koyama, Refinement of prime geodesic theorem,
\textit{Proc. Japan Acad. Ser A Math. Sci.} \textbf{92} (2016),  no. 7, 77--81.

\bibitem{LS:1995} W. Luo and P. Sarnak, Quantum ergodicity of
eigenfunctions on $PSL_{2}(Z)\backslash H^{2}$, \textit{Inst. Hautes \'{E}%
tudes Sci. Publ. Math.} no. 81 (1995), 207--237.

\bibitem{R:1977} B. Randol, On the asymptotic distribution of closed geodesics on compact Riemann surfaces, \textit{Trans. Amer. Math. Soc.} \textbf{233} (1977), 241--247.

\bibitem{SY:2013} K. Soundararajan and M. P. Young, The prime geodesic theorem, \textit{J. Reine Angew. Math.} \textbf{676} (2013), 105--120.
\end{thebibliography}
\end{document}